\documentclass[a4paper,10pt]{article}

\def \Z {{\mathbf {Z}}}

\def \N {{\mathbf {N}}}

\def\u{\bigsqcup}
\def\eps{\varepsilon}
\title{ On Mixing of  Staircase Transformations}
\author{V.V. Ryzhikov}
\date{ Lecture Notes}
\begin{document}
\large
\maketitle

{\bf Staircase  transformations} form a  special subclass of rank one constructions.  Generally a construction is defined by an integer $h_1$, a cut sequence $r_j$    and 
  a sequence $\bar s_j$ of spacer vectors
$$ \bar s_j=(s_j(1), s_j(2),\dots, s_j(r_j-1),s_j(r_j)).$$  Let    $s_j(i)=i-1$,  the corresponding transformations are called staircase.
  Now we   give a more detailed  definition.
Let  a map $T$    be defined (except for the interval  $T^{h_j}E_j$) on $j$-tower 
$$E_j, TE_j, T^2E_j,\dots, T^{h_j}E_j,$$
so $T$ shifts simply  intervals preserving the Lebesgue  measure $\mu$.

 \bf Cutting \rm $E_j$ into $r_j$ subintervals  of the same measure
$E_j=E_j^1\u E_j^2\u  E_j^3\u\dots\u E_j^{r_j},$  
 \bf adding \rm $i-1$ spacers over $i$-th column $E_j^i, TE_j^i ,\dots, T^{h_j}E_j^i,$ for all $i=1,2,\dots, r=r_j$ we  obtain  a partition  
\vspace{2mm}

$E_j^1, TE_j^1, T^2 E_j^1,\dots, T^{h_j-1}E_j^1,T^{h_j}E_j^1,$

$E_j^2, TE_j^2 ,T^2 E_j^2,\dots, T^{h_j-1}E_j^2,T^{h_j}E_j^2,T^{h_j+1}E_j^2,$

$E_j^3, TE_j^3 ,T^2 E_j^3,\dots,T^{h_j-1}E_j^3, T^{h_j}E_j^3,T^{h_j+1}E_j^3,T^{h_j+2}E_j^3,$

$\dots\ \ \dots \ \ \dots$

$E_j^r, TE_j^r ,T^2 E_j^r,\dots,T^{h_j-1}E_j^r, T^{h_j}E_j^r,T^{h_j+1}E_j^r,T^{h_j+2}E_j^r,\dots, T^{h_j+r-1}E_j^r.$
\vspace{2mm}

 \bf Stacking \rm $i+1$-th column over $i$-th column by setting 
$T^{h_j+i}E_j^i = E_j^{i+1}$  for all  $i<r_j$,  
 we get   ($j+1$)-tower 
$$E_{j+1}, TE_{j+1} T^2 E_{j+1},\dots, T^{h_{j+1}}E_{j+1},$$
where 
 $$E_{j+1}= E^1_j, \ \
h_{j+1}+1=(h_j+1)r_j +\sum_{i=1}^{r_j-1}i.$$
Thus,   a measure-preserving map $T$ is defined on all $j$-towers, hence, on their  union $X$ as well.

Adams has proved [1]  the mixing of staircase
transformations  satisfied  the conditions $r_j\to\infty$ and $r_j^2/h_j\to 0$.   

  The author has  presented   (at Paris 6, LPMA,   summer 2000)  his proof of \bf 

 mixing  for all staircases with $r_j\to\infty$ in case of  $\mu(X)=1$. \rm  We reproduce it below.

  Let   a sequence 
   $m_j\in [h_j,h_{j+1})$  satisfy one of the following 
conditions:
$$
 (0)  \ \ \frac{d_j^2}{h_{j}}\to 0; 
\ \ \ \ \ \ \ 
 (C)\ \ \  \frac{d_j^2}{h_{j}}\to c > 0,; 
 \ \ \ \ \ \ \ 
 (\infty) \ \ \ \frac{d_j^2}{h_{j}}\to \infty, 
$$
where  $d_j$ is defined  from the representation 
$$m_j =\sum_{i=1}^{d_j} ( h_{j}+i-1) + t_j, \ \ t_j< d_j+ h_{j}$$ 
(the interval $E^1_{j}$ visited $d_j$ times the  roof of $j$-tower; finally  $T^{m_j}$ maps  it  into $T^{t_j}E_j$; we see that 
$m_j\approx   d_jh_{j}$ for large $d_j$).
Let us fix  sets  $A,B\subset X$  of a positive measure.
We have to prove
$$\mu(A\cap T^{m_j}B) \to  \mu(A)\mu(B), \ \ \ \ \ (T^{m_j} \to_w \Theta)$$ 
(the orthogonal projection into the constants is  denoted by $\Theta$).
We can see that  by a choice of subsequences we come  to one of the cases (0),(C), ($\infty$).
 
The case (0) is exactly  the Adams case.  To prove the mixing in case of (C) we need only to modify Adams' proof.  
 In case of $(\infty)$  we deduce the mixing from the weakly mixing property using ``a progression  of delays''. 
 Now  $s_jr_j>>h_j$, and 
a "geometrical"   picture is the following: a peace-wise linear image  $T^{m_j}E_j$ (or a part of it)  intersects many times the roof of our
$j$-tower and divides the tower  into  certain  domains $D_p$ with   a constant "delay" $d_j(p)$ within $D_p$.
 What is a delay?  If $T^m E^i_j$ is situated in  $T^tE_j$, then $T^m E^{i+1}_j$ will be  in  $T^{t-d}E_j$. So a delay of 
$E^{i+1}_j$ with respect to $E^{i}_j$ is $d=d(i,j)=d_j(p(i))$. (Maybe  it will be  better to say ``upper relative  delay'' of 
$T^mE^{i+1}_j$ with respect to $T^mE^{i}_j$'' is $d(i,j)$.)\vspace{5mm} 
\\ 
{\large
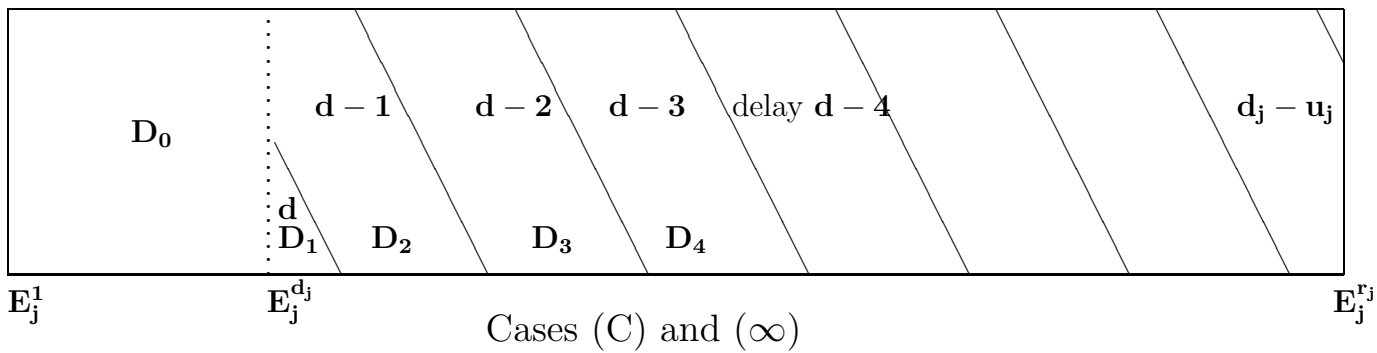
\begin{picture}(0, 100)
\put(0,0){  \line(5, 0){500}}
\put(0,100){  \line(5, 0){500}}
\put(0,0){  \line(0, 1){100}}
\put(500,0){  \line(0,1){100}}

\multiput(100,0)(0,5){20}%
{${\bf .}$}
\put(100,50){  \line(1,-2){25}}
\put(130,100){  \line(1,-2){50}}
\put(190,100){  \line(1,-2){50}}
\put(250,100){  \line(1,-2){50}}
\put(310,100){  \line(1,-2){50}}
\put(490,100){  \line(1,-2){10}}

\put(370,100){  \line(1,-2){50}}
\put(430,100){  \line(1,-2){50}}

\put(50,50){${\bf \Large D_0}$}
\put(105,10){${\bf D_1}$}
\put(105,21){${\bf d}$}
\put(140,10){${\bf D_2}$}
\put(200,10){${\bf D_3}$}
\put(250,10){${\bf D_4}$}
\put(275,60){delay  ${\bf d-4}$}
\put(175,60){  ${\bf d-2}$}
\put(115,60){  ${\bf d-1}$}\put(225,60){  ${\bf d-3}$}
\put(460,60){ ${\bf d_j-u_j}$}
\put(0,-14){ ${\bf E^1_j}$}
\put(97,-14){ ${\bf E^{d_j}_j}$}
\put(496,-14){ ${\bf E^{r_j}_j}$}
\end{picture}
}
{\Large
\begin{center}
Cases (C) and ($\infty$)
\end {center}
}
\vspace{5mm}

 \bf A progression  of delays \rm  is a key property implying the mixing in case of ($\infty$).  Indeed,  there is a standard ``rank one'' approximation technique (see below)  showing
$T^{m_j}\approx_w  \sum_i G_iQ_i,$ 
where  $Q_i= \frac{1}{L_j}\sum_{l=0}^{L_j-1} T^{ld(i,j))}$, and  $\sum_i \|G_i\|<Const$ ( operators $G_i$ are products $\hat Y_i T^{k_i}$, where $\hat Y_i$ are   the multiplications  by  indicators of certain sets $Y_i$).

For \bf a majority of $d(i,j)$ \rm one has 
$Q_i\approx \Theta$,  this implies $T^{m_j}\approx_w  \Theta$,
i.e. the mixing.


{\large
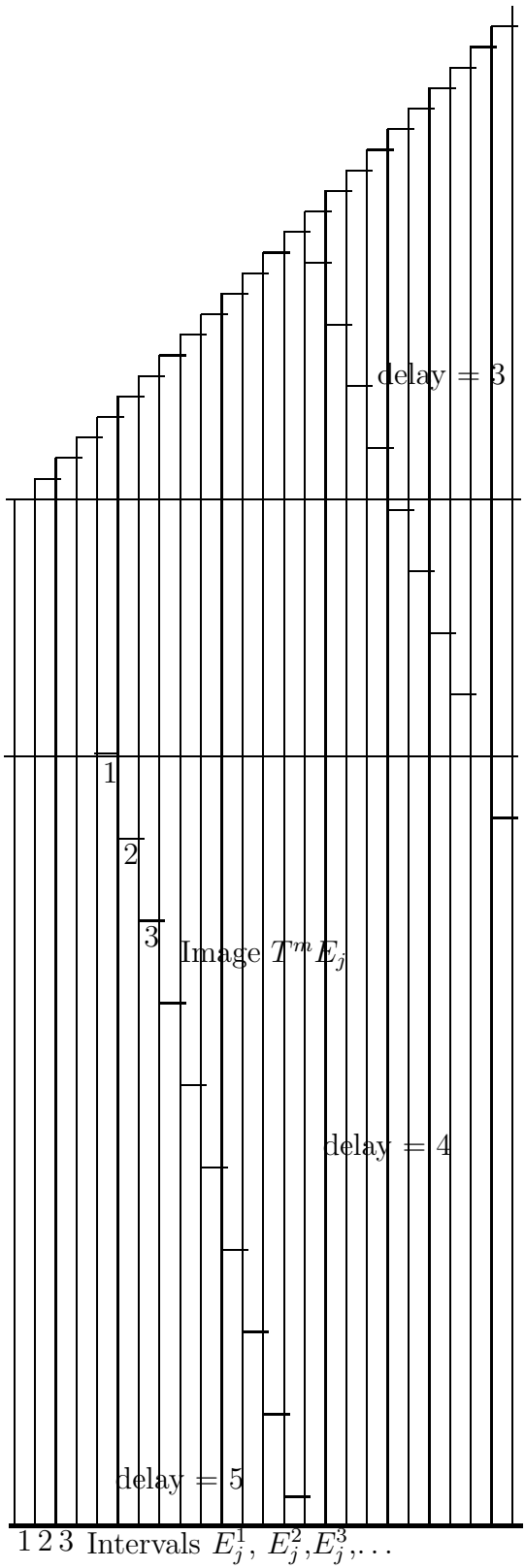
\begin{picture}(650, 560) 
\multiput(0,0)(8,8){25}%
{  \line(0, 1){400}   }
\multiput(0,400)(8,8){24}%
{  \line(1, 0){10}   }
\multiput(0,0)(8,0){25}%
{  \line(0, 1){200}   }
\put(2,0){\line(1,0){198} }
\put(2,1){\line(1,0){198} }
\put(5,-9){1}
\put(13,-9){2}\put(21,-9){3}
\multiput(32,300)(8,-32){10}%
{  \line(1, 0){10}   }
\multiput(112,492)(8,-24){10}%
{  \line(1, 0){10}   }
\put(38,290){1}
\put(46,258){2}\put(54,226){3}
{\large 
\put(228,140){What is a delay?  If $T^m E^i_j$ is situated in  $T^tE_j$, then }
\put(228,125){$T^m E^{i+1}_j$ will be  in  $T^{t-d}E_j$. }
\put(224,110){ This $d$ is an ``upper'' relative  delay}
\put(68,220){Image $T^mE_j$}
\put(220,295){$T^{t_j}E_j$}
\put(220,395){$T^{h_j}E_j$}
\put(220,2){$E_j$}
\put(28,-10){ Intervals $E_j^1$, $E_j^2$,$E_j^3$,$\dots$}
\put(39,15){ delay = 5 }\put(119,145){ delay = 4 }\put(140,445){ delay = 3 }
\put(250,220){$m =\sum_{i=1}^{d_j} ( h_{j}+i-1) + t_j$
  } \put(0,300){\line(1,0){198} }   \put(1,400){\line(1,0){198} } \put(35,301){\line(1,0){8} }  
}
\end{picture}  }

\newpage
{\bf  The case  ($\infty$).  Proof.}  
We define a collection of sets (now  $i$ is not connected with old $i$-s in $E^i_j$)
$$Y_i=\u_{l\in [0,L-1],\\  n\in [0,v]} T^n E^{s(i)+l}_j, \ \
X_i=\u_{ n\in [0,v]} T^n E_j.$$
$$Y_i^+=\u_{l\in [0,L-1],\\  n\in [v',h_j]} T^n E^{s(i)+l}_j, \ \
X_i^+=\u_{  n\in [v',h_j]} T^n E_j.$$

{\large
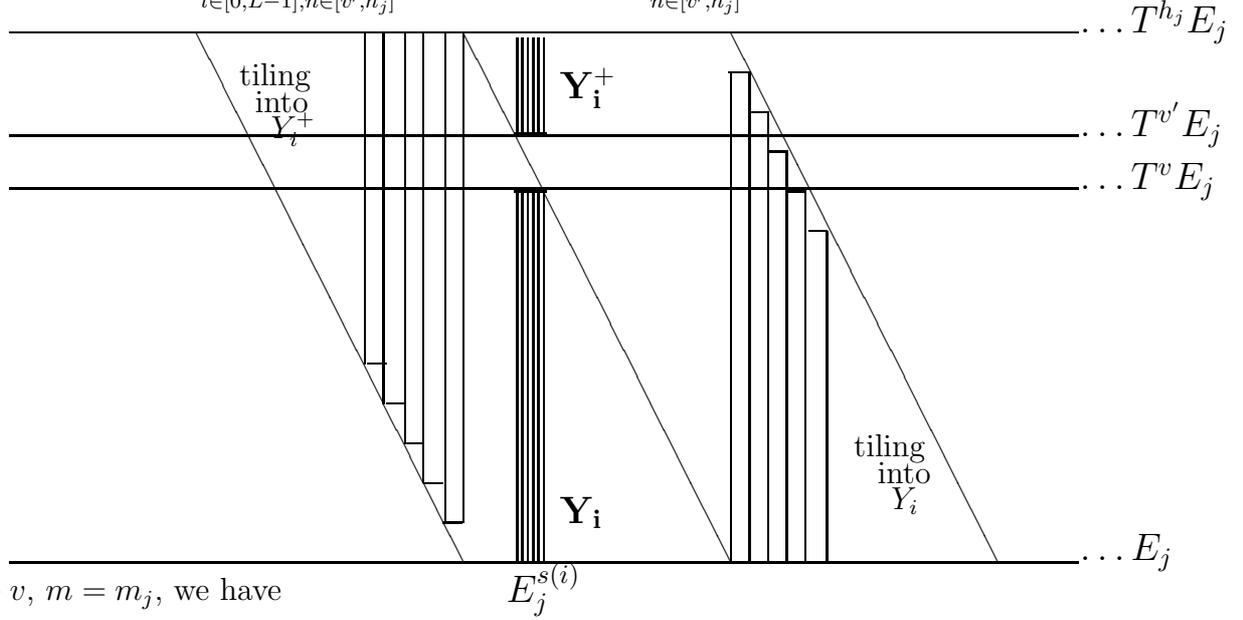
\begin{picture}(0, 200)
\put(50,0){  \line(1, 0){400}}{\Large \put(450,1){  $\dots E_j$}}
\put(50,200){  \line(5, 0){400}}{\Large \put(450,200){  $\dots T^{h_j}E_j$}}

\put(120,200){  \line(1,-2){100}}
\put(220,200){  \line(1,-2){100}}

\put(220,200){  \line(0,-1){185}}
\put(223,15){\line(-1,0){7}}
\put(213,200){  \line(0,-1){185}}

\put(216,30){\line(-1,0){7}}
\put(205,200){  \line(0,-1){170}}
\put(209,45){\line(-1,0){7}}
\put(198,200){  \line(0,-1){155}}
\put(202,60){\line(-1,0){7}}
\put(190,200){  \line(0,-1){140}}

\put(195,75){\line(-1,0){7}}
\put(183,200){  \line(0,-1){125}}
\put(140,180){tiling }
\put(142,170){ into }
\put(148,160){ $Y^+_i$}

\put(320,0){  \line(0,1){185}}
\put(323,185){\line(1,0){7}}

\put(327,0){  \line(0,1){185}}
\put(331,170){\line(1,0){7}}
\put(334,0){  \line(0,1){170}}
\put(338,155){\line(1,0){7}}
\put(341,0){  \line(0,1){155}}
\put(345,140){\line(1,0){7}}
\put(348,0){  \line(0,1){140}}
\put(353,125){\line(1,0){7}}
\put(356,0){  \line(0,1){125}}

\put(370, 40){tiling }
\put(375,30){ into }
\put(380,20){ $Y_i$}

\multiput(250,0)(-2,0){6}%
{  \line(0,1){140}}
\put(239,162){  \line(1,0){12}}
\put(239,140){  \line(1,0){12}}
\put(50,141){  \line(1,0){400}} {\Large \put(450,141){  $\dots T^vE_j$}}
\put(50,161){  \line(1,0){400}} {\Large \put(450,161){  $\dots T^{v'}E_j$}}
\multiput(250,162)(-2,0){6}%
{  \line(0,1){36}}
\put(320,200){  \line(1,-2){100}}
{\Large
\put(260,15){${\bf Y_i}$}
\put(260,175){${\bf Y_i^+}$}
\put(240,-15){${ E_j^{s(i)}}$}
}
\end{picture}
}
\\  Let $k=h_j-v$,  $m=m_j$,  we have
$$\mu(T^m A\cap B|Y_i)= \frac{1}{L}\sum_{l=0}^{L-1}\mu(T^kT^{ld}A\cap B|Y_i)= $$
$$=\frac{1}{L}\sum_{l=0}^{L-1}\mu(T^{ld}A\cap T^{-k}B|T^{-k}X_i)= $$
$$\mu(A)\mu(B) - \int \left(\frac{1}{L}\sum_{l=0}^{L-1}T^{ld}f
\right)T^{-k}(\hat B\hat X_i) d\mu\,/\mu(X_i),$$
where $f= \chi_{{}_ B} -\mu(B)$.  

We fix a very small $\delta>0$  and consider from now only     $v$  satisfied $\delta h_j < v < (1-\delta)h_j$.  Since $\mu(X_i)>\delta$ from now, we get 
$$|\mu(T^m A\cap B\cap Y_i) - \mu(A)\mu(B)\mu(Y_i)|\leq \mu(Y_i)\delta^{-1}\left|\left| \frac{1}{L}\sum_{l=0}^{L-1}T^{ld}f\right|\right|.\eqno (1)$$

 To prove the mixing in the case ($\infty$)  let's recall the Blum-Hanson inequality:
if $\langle T^{ld}f, f\rangle <\eps$ for all $l=1,2,\dots , L$, then 
$$ \left|\left|\frac{1}{L}\sum_{l=0}^{L-1} T^{ld}f\right|\right|  ^2 \leq \frac{1}{L^2}(\eps L(L-1) +L\|f\|^2)< \eps+\frac{1}{L}.$$

\bf Lemma 1. \it  Let  $T$ be weakly mixing,  $\eps>0$  and $f\in L_2^0$, then    
$$\frac {|\{d\in I_j=[d_j, d_j+u_j]\ :\ \langle T^{id}f, f\rangle <\eps,  \ i=1,2,\dots , L-1\}|}{u_j}\ \to \ 1$$  
as $u_j\to\infty$. (We say that a majority of $d$  consists from  $\eps$-good $d$-s.)
\rm

Indeed,  we know that 
$|I_j|^{-1}\sum_{d\in I_j} S^d \to_w \Theta$ for ergodic $S$.  Let $f$ be a real function. 
If for all  $d\in J_j\subset I_j$ we have $\langle S^{d}f, f\rangle \geq \eps>0,$ 
(or  $\langle S^{d}f, f\rangle \leq -\eps$ ), 
then $$|J_{j'}|^{-1}\sum_{d\in I_j} S^d \to_w P\neq\Theta.$$
Let $|J_{j'}|/|I_{j'}|\to a$.   We get $\Theta =aP+(1-a)P'$, where the Markov operators  $P, P'$ commute with $S$.
Since $S$ is weakly mixing,  $P\neq \Theta$ implies $a=0$  ($\Theta$ is an extreme point in Markov centralizer of $S$).  
 So, we showed  that for any fixed  $i$  most of $d$ satisfied  
$\langle T^{id}f, f\rangle <\eps$.  Setting $S=T^l$ we prove lemma.

Lemma 1  says that a majority of  delays $d\in [d_j, d_j+u_j]$ consists from  
 $\eps$-good ones, hence,    from  (1)   we  get   
$$\mu(T^m A\cap B\cap Y_i) = \mu(A)\mu(B)\mu(Y_i) + \mu(Y_i) \delta^{-1}\sqrt{\eps+\frac{1}{L}} $$
for all $Y_i$ ($Y^+_j$ as well) from  all domains $D_p$ with $\eps$-good delays (a measure of all other domains is vanishing !). 
Thus, Applying Lemma 1 we  cover a majority  of $D_1, D_2, D_3,\dots$  by a union of  sets $Y_i$,  $Y^+_i$ with $\eps$-good delays $d=d_j(i)$.  Finally let $L\to\infty$, $\eps, \delta\to 0$ very slowly.  
 The mixing in case of ($\infty$) is proved.

{\bf  The case  (C).  Proof.} Now we have  $d_jr_j\approx C h_j$.  From
$$ h_{j-1}^2 >>Ch_j \approx  d_jr_j> d_j^2$$
we get  
$$ h_{j-1} >> d_j.$$
So, there is $q$ such that $qd_j\in [h_{j-1}, 2h_{j-1}]$.  The intervals $[h_{j-1}, Lh_{j-1}]$ are (asymptotically) mixing. 
For $d\in [d_j -1, d_j + N]$, where some  $N>C$ is fixed, we  
have 
 $$T^{qd}, T^{2qd},\dots, T^{Lqd}\approx_\eps \Theta$$
that means: $\langle T^{qd}f, f\rangle<\eps$  for all $ i=1,2,\dots , L$;   $f= \chi_{{}_ B} -\mu(B)$. 

Defining again the 
 ``rakes'' 
$$Y_i=\u_{l\in [0,L-1],\\  n\in [0,V]} T^n E^{s(i)+lqd}_j, \ \
Y_i^+=\u_{l\in [0,L-1],\\  n\in [v',h_j]} T^n E^{s(i)+lqd}_j,  $$
 using  
 the Blum-Hanson-Adams trick we get the mixing on our rakes:
$$\mu(T^m A\cap B|Y_i)= \frac{1}{L}\sum_{l=0}^{L-1}\mu(T^kT^{lqd}A\cap B|Y_i)\approx \mu(A)\mu(B). \eqno BHA-trick$$
We have (trivially) the mixing on $D_0$ ( amusing,  we can  use   rakes  for   $D_0$ with a delay $d= 1$).  We cover $X\setminus D_0$ by rakes,  we obtain the mixing. Thus,
$\mu(T^{m_j} A\cap B)\approx \mu(A)\mu(B)$
is proved in case of (C).

{\bf  The Adams case  (0).  Proof.}  Again we define  $Y_i, Y_i^+$ as above.
\vspace{5mm}
\\
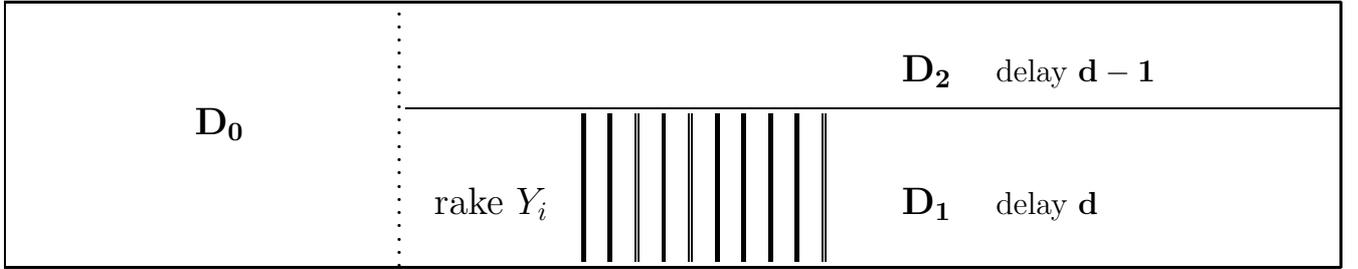
\begin{picture}(0, 100)
\put(0,0){  \line(5, 0){500}}
\put(0,100){  \line(5, 0){500}}
\put(0,0){  \line(0, 1){100}}
\put(500,0){  \line(0,1){100}}

\multiput(220,2)(10,0){10}%
{\line(0,1){56}}
\multiput(221,2)(10,0){10}%
{\line(0,1){56}}

\multiput(150,0)(0,5){20}%
{${\bf .}$}
\put(150,60){  \line(1,-0){350}}
{\Large
\put(75,50){${\bf  D_0}$}
\put(335,70){  ${\bf D_2}$}
\put(160,20){ rake  $Y_i$}
\put(335,20){  ${\bf D_1}$}}
\put(375,70){delay  ${\bf d-1}$}
\put(375,20){delay  ${\bf d}$}
\end{picture}
\vspace{1mm}

\vspace{1mm}
Now  we find $q$  such that    $h_{j_p-1}<qd_j<2h_{j_p-1}$ for some $p=p_j$ (if $d_j\to\infty$, then $p_j  \to \infty$; if $d_j=d$ is constant, we use simply the ergodicity of $T^d$).   Adams  showed the existence of  such
$p$ as follows.  

If $ h_p \approx Cr_jd_j$, then  for $h_{p-1}$ we have   $h_{p-1}^2>>r_jd_j$, so  $r_jd_j>>h_{p-1}>>d_j$, and we find a desired $q$ for $h_{p-1}$.       

If   $h_{p} <<r_jd_j \leq h_{p+1}$, then  $r_jd_j  <<h_{p}^2$,  so  $d_j<<h_{p} <<r_jd_j$.  We find $q$, consider mixing intervals 
$[h_p, Lh_p]$ and 
 prove  mixing on $Y_i$  via  BHA-trick.  Then  covering  $D_1$ and $D_2$ by sets $Y_i,Y_i^+$ provides the mixing.
\vspace{2mm}

{\normalsize      
We note that Adams' method has been   generalized in several directions  by D.Creutz, C.Silva, and by the present author.  We shall discuss some generalizations later.    D.Creutz and  C.Silva  have got their proof of the mixing for staircase transformations (see \it  Mixing on rank-one transformations, \rm Studia Math. 199, 2010, 43-72).  Our  method, we think,  is more  ``direct'' and more close to Adam's proof.
\vspace{2mm}}

\bf On "Infinite" Transformations. \rm 
Let $\mu(X)=\infty$, i.e. we  consider a situation in case of  $\sum_j r_j/h_j=\infty$.
The above proof is  quite suitable for  the infinite measure case.
\vspace{2mm}
\\
\bf THEOREM. \it If $r_j/h_j\to 0$ and $r_j\to\infty$, then the corresponding staircase transformation is mixing:
for any  $A,B$ of a finite measure $\mu(T^mA\cap B)\to 0$ as $m\to\infty$.\rm

Proof.  We repeat the above proof only replacing  $\Theta$ by $0$,  and apply  the following simple lemma instead of Lemma 1. 
 
\bf Lemma 2. ($\eps$-good $d$-s prevail) \it  If $T$ is ergodic, then  for any 
 $L$, $\eps>0$,  $f\in L_2(X)$,  $u_j\to\infty$ we have
$${u_j^{-1}} {\left|\{d\in [d_j, d_j+u_j]\, :\, \langle T^{ld}f, f\rangle <\eps,  \ l=1,2,\dots , L-1\}\right|}\ \to \ 1.$$ 
\rm

Let's pay finally attention to one of open problems.
\\ 
\bf  Problem. \it    Prove the mixing for the "$r_j=h_j$" staircase transformation.\rm
\\
\bf  Conjecture. \it Any staircase transformation   is mixing as  $r_j\to\infty$.\rm \vspace{2mm}
\\ \normalsize
\bf [1] \rm T.M. Adams. \it  Smorodinsky's conjecture on rank one systems, \rm Proc. Amer. Math. Soc.
 {\bf 126} (1998), 739-744.
\\ \\
vryzh@mail.ru
\end{document}